%% file: crys-lex.tex
\documentclass{amsart}
\include{headerstemplate}

\begin{document}
\title{Mochizuki's Crys-Stable Bundles: A Lexicon and Applications}
\author{Brian Osserman}
\begin{abstract}
Mochizuki's work on torally crys-stable bundles \cite{mo3}
has extensive implications for the
theory of logarithmic connections on vector bundles of rank $2$ on curves,
once the language is translated appropriately. We describe how to carry out 
this translation, and give two classes of applications: first, one can
conclude almost immediately certain results classifying Frobenius-unstable 
vector bundles on curves; and second, it follows from the results of
\cite{os10} that one also obtains results on rational functions with 
prescribed ramification in positive characteristic. 
\end{abstract}
\thanks{This paper was supported by a fellowship from the Japan Society for
the Promotion of Science.}
\maketitle

\section{Introduction}

Mochizuki's theory of torally crys-stable bundles and torally indigenous
bundles developed in \cite{mo3} has, after appropriate translation,
immediate implications for logarithmic connections on vector bundles of rank
2 on curves. This in turn has immediate implications to a subject which has
recently been studied by a number of different people (see, for instance, 
\cite{l-p}, \cite{l-p3}, \cite{j-r-x-y}, \cite{os11}): Frobenius-unstable
vector bundles, and by extension the generalized Verschiebung rational map 
induced on moduli spaces of vector bundles by pulling back under Frobenius.
Furthermore, together with the results of \cite{os10}, one can use
Mochizuki's work to describe rational functions with prescribed ramification
in positive characteristic. This relationship provided the original 
motivation for the
ultimately self-contained arguments of \cite{os7}, but also goes further by
giving a finiteness result not yet known by direct arguments.

However, because the theory of \cite{mo3} was aimed towards unrelated
$p$-adic applications, it does not contain any translation of its results
back into the language of vector bundles with connection or
Frobenius-unstable vector bundles, and until recently remained outside the
circle of literature treating these topics. The aim of this paper is 
therefore primarily expository: we will give a lexicon of Mochizuki's 
language and a survey of certain key results, and then detail the 
applications described above. The arguments in Section \ref{s-proofs} were suggested almost uniformally by Mochizuki, but any blame for the particular translation in Section \ref{s-lexicon} of Mochizuki's work into the language of vector bundles lies solely with the author.

\subsection*{Acknowledgements}

I would like to thank Shinichi Mochizuki for bringing the relevant results of
\cite{mo3} to my attention, and assisting me greatly in translating them
into more familiar language. I would also like to thank Max Lieblich for his helpful conversations.

\section{The Lexicon}\label{s-lexicon}

Fix $g,r \geq 0$ with $2g-2+r>0$, and $p$ an odd prime. Let $S$ be a local
base scheme, defined over the residue field $k$ of its closed point. We 
refer to a curve $C$ over $S$ or bundle on $C$ as ``constant'' if it is
obtained by pullback from a $C_0$ or bundle on $C_0$ over $\Spec k$. Now,
let $C$ be a smooth 
proper curve of genus $g$, with $r$ marked points $P_1, \dots, P_r$. 
For technical reasons, we require only that $S$ is Noetherian and strictly 
Henselian of
characteristic $p$, and $C$ constant. However, for conceptual purposes, it 
suffices to take $S = \Spec k$, for an algebraically closed field $k$ of
characteristic $p$. We begin by setting terminological conventions for some standard concepts.

\begin{defn} Given a vector bundle $\E$ on $C$, a {\bf logarithmic}
$S$-connection on $\E$ is an $\O_S$-linear map $\nabla:\E \rightarrow 
\Omega^1_{C/S}(\sum_i [P_i]) \otimes \E$, satisfying the connection rule
$\nabla(fs) = f \nabla(s) + df \otimes s$. For any $i$, a logarithmic 
$S$-connection has a canonical {\bf residue} $\Res_{P_i} \nabla
\in \End(\E|_x)$ induced by the residue map $\Omega^1_{C/S}(P_i) 
\rightarrow \Gamma(S, \O_S)$.
\end{defn}

\begin{defn} Given a vector bundle $\E$ on $C$, a sub-bundle $\F$ of $\E$ is said to be {\bf destabilizing} (respectively, {\bf weakly destabilizing}) if the slope of $\F$ is strictly greater than (resp., greater than or equal to) the slope of $\E$.
\end{defn}

Because $S$ has characteristic $p$, there is a notion of $p$-curvature 
associated to $\nabla$ (see \cite[\S 5]{ka1}), which allows us to define:

\begin{defn} We say that a logarithmic $S$-connection has {\bf $p$-trivial
determinant} if the determinant connection $\det \nabla$ on $\det \E$ is a
regular connection with vanishing $p$-curvature.
\end{defn}

As the name suggests, this definition is intended to generalize trivial
determinant to cases where the stricter notion wouldn't make sense -- for
instance, when $\E$ has determinant of degree $p$. We also remark that the 
regularity of the determinant connection implies that the residues of 
$\nabla$ have vanishing trace.

For each $i$, fix a $\rho_i \in \Gamma(S, \O_S)/\pm 1$.
The fundamental object we will deal with is, in Mochizuki's language, the
{\bf torally crys-stable bundle} (see \cite[Def. I.1.2, p. 89]{mo3}). In our
lexicon, we have the following equivalence:

\begin{prop}\label{crys-stab} A torally crys-stable bundle on $C/S$ of radii
$\{\rho_i\}$ is
equivalent to an equivalence class of vector bundles $\E$ of rank $2$ on 
$C$, together with logarithmic $S$-connection $\nabla$ having $p$-trivial
determinant, under the 
equivalence obtained by allowing tensoring by any line bundle with 
$S$-connection. The pair $(\E, \nabla)$ must also satisfy the
following additional conditions, which are invariant under our equivalence
relation:
\begin{ilist}
\itm for each $i$, $\Tr ((\Res_{P_i} \nabla)^2)= 2\rho_i ^2$; 
\itm at the closed point of $S$, $\Res_{P_i} \nabla$ is non-zero for all
$i$;
\itm after restriction to an arbitrary geometric point of $S$, there is no 
weakly destabilizing sub-bundle of $\End^0(\E)$ which is horizontal with respect to the connection induced by $\nabla$.
\end{ilist}
\end{prop}

Since $\rho_i$ was only chosen up to sign, we see that the radii of a
torally crys-stable bundle are well-defined. We also note that if for any
$i$,
$\Res_{P_i} \nabla$ is diagonalizable, condition (i) is vacuous, and the
radius will be the corresponding eigenvalues (which will necessarily be
paired with opposite sign). Finally, we remark that although condition (ii)
is very closely related to the non-existence of a horizontal destabilizing
sub-bundle of $\E$, stating the relationship precisely is complicated, and
as stated below, will not be necessary.

We now give equivalent definitions of certain other terms of Mochizuki's
theory which will be relevant. First, we define:

\begin{defn} Let $(\E, \nabla)$ be a logarithmic $S$-connection on $C$, and
$\L$ a line sub-bundle of $\E$. Then the {\bf Kodaira-Spencer map} 
$$\kappa: \L \rightarrow \E/\L \otimes \Omega^1_{C/S}
(\sum_i [P_i])$$
associated to $(\E, \nabla)$ and $\L$ is obtained by applying $\nabla$ to
the natural inclusion and quotient maps. When $\E$ is unstable, we will
refer to the Kodaira-Spencer map of $(\E, \nabla)$ to mean the
Kodaira-Spencer map associated to $(\E, \nabla)$ and the destabilizing
sub-bundle of $\E$.
\end{defn}

The second definition makes sense by Lemma \ref{destab-unique}.
Note that a line sub-bundle $\L$ being horizontal for
$\nabla$ is equivalent to the vanishing of its Kodaira-Spencer map.

\begin{prop}\label{level} Let $(\E, \nabla)$ be a representative of a
torally crys-stable bundle. Then $(\E, \nabla)$ has {\bf level} $\ell$ if:
\begin{ilist}
\itm in the case $\ell=0$, the restriction of $\E$ to every geometric point
of $S$ is semistable;
\itm in the case $\ell>0$, $\E$ has a line sub-bundle $\L$ of degree
$\ell+\frac{1}{2}\deg \E$ such that the Kodaira-Spencer map associated to
$(\E, \nabla)$ is non-zero, and is an isomorphism at the $P_i$.
\end{ilist}

Furthermore, the present condition (ii) supercedes condition (iii) of
Proposition \ref{crys-stab}.
\end{prop}

Thus, the level of a torally crys-stable bundle is independent of $\nabla$,
and simply measures the failure of $\E$ to be semi-stable. Note however that
unless $S = \Spec k$, it is not automatically the case that any torally 
crys-stable bundle has a level. 

It immediately follows (see \cite[Lem. I.3.4, p. 104]{mo3}) that we can 
define:

\begin{defn}\label{indig} A {\bf torally indigenous bundle} is a torally
crys-stable bundle of level $\frac{1}{2}(2g-2+r)$.
\end{defn}

Finally, we can define:

\begin{prop}\label{dormant} A torally crys-stable bundle represented by
$(\E, \nabla)$ is {\bf dormant} if the $p$-curvature of $\nabla$ is zero.
\end{prop}

Note that here we use the hypothesis of Proposition \ref{crys-stab} that our connections have $p$-trivial determinant, as we otherwise would not have that vanishing $p$-curvature is preserved under equivalence. 

The results of Mochizuki central for our purposes can now be stated:

\begin{thm}\label{moch1} (Mochizuki) The stack of dormant
torally indigenous bundles is finite flat over $\overline{\M}_{g,r}$, and
\'etale over points corresponding to totally degenerate curves.
\end{thm}

\begin{thm}\label{moch2} (Mochizuki) Suppose that $p \geq 2g+r-2\ell$, and
$\ell>0$. Then the stack of dormant torally crys-stable bundles of level 
$\ell$ is flat over $\overline{\M}_{g,r}$ of relative dimension 
$2g-2+r-2\ell$.
\end{thm}

\section{Proofs}\label{s-proofs}

We begin with a standard lemma, observing that it holds over an arbitrary
base:

\begin{lem}\label{destab-unique} Let $C/S$ be a smooth, proper curve
over a connected scheme $S$, and $\E$ a vector bundle of rank
$2$ on $C$. Suppose that $\L \subset \E$ is a weakly destabilizing line 
sub-bundle of $\E$. Then $\deg \L$ is uniquely determined,
and if $\L$ is in fact destabilizing, $\L$ itself is uniquely determined.
Furthermore, the same assertions holds for sub-bundles of $\End^0(\E)$ of a given rank.
\end{lem}

\begin{proof} The only observation required is that even in the general
relative setting, Nakayama's lemma gives the statement that a line bundle 
of negative degree has no non-zero global sections. The assertions then 
follow by standard arguments involving composing the inclusion of one line 
sub-bundle with the quotient map induced by another, and vice versa.
The argument in the case of destabilizing line sub-bundles of $\End^0(\E)$
is similar, and given in \cite[Proof of Lem. I.3.5, p. 105]{mo3}. In fact, one checks that this argument yields the stronger statement, if in its notation we let $\L_2$ be a destabilizing line sub-bundle, and $\L_1$ a weakly destabilizing line sub-bundle. Finally, the rank $2$ case follows from the rank $1$ case by self-duality.
\end{proof}

Next, we show that our definition of torally crys-stable agrees with
Mochizuki's, at least over a strictly Henselian base. We begin with some 
general remarks on projective bundles and connections.

First, given a projective bundle $P$ on a relative smooth curve $C/S$, it is well-known that \'etale locally on $S$, we may write $P=\P(\E)$ for some 
vector bundle $\E$ on $C$. Indeed, this is equivalent to finding a line bundle on $P$ which is $\O(1)$ on fibers, from which one can obtain $\E$ via push-forward. By Tsen's theorem one can remove the Brauer obstruction for such an $\E$ on any special fiber after \'etale base change, one can extend to the complete local ring on the base by vanishing of obstructions in deformation theory, and finally apply the Artin approximation theorem to construct $\E$ after a further \'etale base change. Given this construction, since $\E$ is unique up to tensoring by line bundles, and by the invariance of degree in families, we may define:

\begin{defn} Let $P$ be a $\P^{r-1}$-bundle on a smooth relative curve 
$C/S$. Then the {\bf degree class} of $P$ is defined to be the congruence class modulo $r$ of the degree of any $\E$ on $C$ with $P\cong \P(\E)$.
\end{defn}

Note that the degree class could be defined more intrinsically as the class induced by $P$ in $H^2_{\et}(C, \mu_r)$; however, the more elementary definition will also be more immediately useful.

There are a number of issues involving translation between the formalities of the Grothendieck point of view on connections and the classical definitions, and similarly for the formalism of log structures. However, such translations are well-known, and we refer the reader to \cite[Prop. 2.9]{b-o}, \cite{ka4}, and \cite[p. 13]{og2}. 
For our arguments, it will be convenient to use the following naive definition. 

\begin{defn} Let $\E$ be a vector bundle on $C/S$. Then $\nabla$ is an {\bf
$S$-connection on $\P(\E)$} if it is a global section of the sheaf $\PConn$
obtained as the sheafification of the presheaf whose sections over any $U
\subset C$ are equivalence classes of $S$-connections on $\E|_U$, defined up
to tensoring by connections on line bundles. A {\bf rational} $S$-connection $\nabla$ on $\P(\E)$ is a section of $\PConn$ over some $U \subset C$, and the {\bf order} of a pole of $\nabla$ at $P \in C \smallsetminus U$ is defined to be the minimum order of the pole at $P$ of all rational connections on $\E$ representing $\nabla$ on $U' \subset U$. If $\nabla$ has only simple poles, it is {\bf logarithmic}.
\end{defn}

Such a connection on $\E$ is easily checked to yield a connection in the Grothendieck sense on $\P(\E)$ as a scheme over $C$, while one may go backwards \'etale locally, at least when $p$ is prime to $r$, by finding an $\O(1)$ on $\P(\E)$ which is preserved by the connection, thereby giving a connection on its push-forward to $C$, which is isomorphic to $\E$.

\begin{prop}\label{deproj} Suppose that $p$ does not divide $r$. Let 
$(P, \nabla^P)$ be a projective bundle of dimension $r-1$ with rational 
connection on $C/S$. Let $\L$ be a line bundle whose degree is the degree 
class of $P$, and with connection $\nabla^0$. Then there exists an $\E$ of
determinant $\L$ on $C/S$ such that $\P(\E) \cong P$, and for any such $\E$,
there is a unique rational connection $\nabla$ whose determinant is
$\nabla^0$, and such that $\nabla$ recovers $\nabla^P$ under
projectivization. Further, the locations and orders of the poles of 
$\nabla$ are the same as those of $\nabla^P$.
\end{prop}

\begin{proof} Let $\E'$ be any vector bundle such that $\P(\E') \cong P$.
Since $\L$ is of the appropriate degree class, and multiplication by $r$ is
\'etale on the Jacobian, we find that there exists an $\L'$ with
$\L'^{\otimes r} \cong \L \otimes (\det \E)^{-1}$, and setting $\E = \E'
\otimes \L'$ will give us $\det \E \cong \L$. Next, since $p$ does not 
divide $r$, $\nabla$ will certainly be uniquely determined by its 
determinant and projectivization,
if it exists. Moreover, everything is compatible with restriction, so by
uniqueness it suffices to produce $\nabla$ on an open cover of $C$, since
the restrictions on overlaps will then have to glue. Let $U$ be an open set
on which $\nabla^P$ is represented by some $\nabla^U$ on $\E|_U$, having the 
same locations and orders of poles as $\nabla^P$. We then see that by
tensoring with $(\O_U, \frac{1}{r}(\nabla^0-\det \nabla^U))$, we obtain the
desired $\nabla$ on $U$. Finally, this clearly cannot increase pole orders, 
but neither can it decrease them by the 
definition of the order of a pole of a projective connection.
\end{proof}

\begin{proof}[Proof of Proposition \ref{crys-stab}] Given
a vector bundle with connection $(\E, \nabla)$, we first want to check that
conditions (i)-(iii) of the proposition are equivalent to conditions (1)-(3)
of \cite[Def. 1.2, p. 89]{mo3} for the projectivization $\P(\E)$ with the 
induced connection. For conditions (i) and (ii), one need only know that
Mochizuki's monodromy operator $\mu_i$ is simply the residue of
$\nabla$ at $P_i$. Then, because $\Ad(\P(\E)) \cong \End^0(\E)$, the
equivalence for condition (iii) is trivial. Now, condition (4) is vacuous 
since we have restricted to smooth curves, so
we find that given a vector bundle with connection $(\E, \nabla)$ satisfying
the condition of the proposition, we get a torally crys-stable bundle simply
by projectivizing. 

Conversely, we want to show that, since our base is
strictly Henselian, given a torally crys-stable bundle $(P, \nabla)$ we can
find a vector bundle with connection having $p$-trivial determinant, and
such that the projectivization recovers $(P, \nabla)$. By Proposition
\ref{deproj}, we need only produce a line bundle $\L$ of the
appropriate degree class, having a connection with vanishing $p$-curvature.
But producing an $\L$ of the appropriate degree class is trivial, and we can
and do require further that it be constant, and that $p$ divide its degree.
Using the surjectivity of the Verschiebung map on
Jacobians over an algebraically closed field and the fact that $\L$ is
constant, we find that $\L \cong F^*
\L'$ for some $\L'$, and the induced canonical connection gives us the
desired connection with vanishing $p$-curvature.
\end{proof}

We next have:

\begin{proof}[Proof of Proposition \ref{level}] We wish to compare our
criteria to \cite[Def. I.3.2, p. 103]{mo3}. In light of the translation
provided by the above argument, and using the equivalence between line
sub-bundles of a rank-$2$ vector bundle and sections of its 
projectivization, one need only check two statements: first, that 
if $h_{\L}$ denotes the section of $\P(\E)$ associated to a line sub-bundle 
$\L$ of $\E$, that $h_{\L}^* \omega_{\P(\E)/C} \cong \Hom(\E/\L,\L)$; and 
second, that condition (ii) of the proposition is
equivalent to condition (2) of the definition together with condition (iii) 
of Proposition \ref{crys-stab}. For the first, we consider the dual 
statement, and note that $h_{\L}^* \tau_{\P(\E)/C}$ is the first-order 
deformation space
of $h_{\L}$, which is naturally identified with $\Hom(\L,\E/\L)$.

For the second, the main observation is the following: given a destabilizing
line bundle $\L \subset \E$, we obtain a destabilizing line sub-bundle $\L^0
\subset \End^0(\E)$ by considering the traceless endomorphisms mapping $\L$
to $0$. Indeed, it is easy to see that $\L^0 \cong \Hom(\E/\L, \L)$, which
has positive degree, and since $\End^0(\E)$ has degree $0$, $\L^0$ is
destabilizing. Further, we note that it is easy to check that $\L$ is
horizontal for $\nabla$ if and only if $\L^0$ is horizontal for the
connection induced by $\nabla$ on $\End^0(\E)$. Now, by Lemma
\ref{destab-unique}, we see that the existence of a non-horizontal
destabilizing $\L$ in condition (ii) of the present proposition precludes
the existence of a horizontal weakly destabilizing line sub-bundle, which by
self-duality of $\End^0(\E)$ precludes the existence of a horizontal
weakly destabilizing bundle of any rank $2$, also. Conversely, condition (iii) of
Proposition \ref{crys-stab} implies that $\L$ cannot be horizontal, so we
obtain the desired equivalence.
\end{proof}

Before discussing our definition of $p$-curvature, we need to discuss the lifting of Grothendieck-perspective connections to higher-order PD-neighborhoods. This is completely standard in the case of classical connections on a sheaf of modules on $C/S$; see, e.g., 
\cite[Thm. 4.8]{b-o}. Now, suppose we have a connection $\nabla$ on some
object $P$ over $C$, which for simplicity we assume to be a separated scheme
over $C$; that is, if $C^{[2]}$ denotes the first-order neighborhood of the
diagonal in $C \times _S C$ (which is the same with or without divided
powers), $\nabla$ is an $\O_{C^{[2]}}$-linear isomorphism $p_1^* P \cong p_2
^* P$, which gives the identity when restricted to the diagonal. We wish to
show that (since $C$ is a curve, so any connection is automatically
integrable), we obtain a unique lifting to $n$th-order neighborhoods of the
diagonal in $C \times^{\PD}_S C$ for any $n$. Let $U \rightarrow P$ be any
affine cover; in particular, $U$ is affine over $C$, so $\O_U$ is
quasi-coherent over $C$, and the connection on $P$ induced a classical
connection on $\O_U$, which furthermore is compatible with the algebra
structure on $\O_U$. Now, by the classical theory, we can lift this
connection uniquely to $n$th-order neighborhoods, with the lift given by the formal expression $e^{(\nabla_{\frac{d}{dt}})\otimes(1 \otimes t - t \otimes 1)}$ for some local coordinate $t$: that is to say, in the notation of \cite{b-o}, we take the map $\E \rightarrow \E \otimes \P^n$ given by $s \mapsto \sum_{i=0}^{n} (\nabla_{\frac{d}{dt}})^i s \otimes (1 \otimes t - t \otimes 1)^{[i]}$, and extend by scalars to obtained the desired $\P^n$-linear isomorphism $\P^n \otimes \E \rightarrow \E \otimes \P^n$. One checks that the product rule for exponentials implies 
that this lift preserves the compatibility with the algebra structure. By 
uniqueness, this 
lifting necessarily agrees on restriction to intersections of open sets in 
the cover, so it follows that one obtains a lifting of the connection on 
$P$ to $n$th-order as well.

We are now ready to define $p$-curvature from the Grothendieck perspective. the idea is to consider the following picture:
$$\xymatrix{{C \times^{\PD}_S C} \ar[d] & {V(\I)} \ar[d]\ar@{_{(}->}[l] & {V(\I,
\J^{[p+1]})} \ar@{_{(}->}[l]\\
{C \times _S C} & {C} \ar@{_{(}->}[l]}$$
where the square is Cartesian, so that $\I$ is by definition the ideal in $C
\times^{\PD}_S C$ generated by the ideal defining the diagonal in $C
\times_S C$. If we denote by $\J$ the ideal of the diagonal in 
$C \times^{\PD}_S C$, we have $\J = \cup_{i \geq 1} \I^{[i]}$. Now, the main 
observation is that $V(\I,\J^{[p+1]})$ has structure sheaf isomorphic to 
$\O_C \oplus F_C^* \omega_{C/S}$, with multiplication determined by setting 
the product of any two elements of $F_C^* \omega_{C/S}$ to $0$. Given 
$\nabla$ an isomorphism between $p_1^* P$ and $p_2^* P$ on the first-order
neighborhood of the diagonal in $C \times_S C$, we can pull back to $C
\times^{\PD} _S C$, to obtain such an isomorphism on $V(\J^{[2]})$, and we 
can then lift to $(p+1)$st order, obtaining an isomorphism on
$V(\J^{[p+1]})$, which we can restrict to $V(\I, \J^{[p+1]})$. On the other
hand, we also have the trivial isomorphism between $p_1^* P$ and $p_2^* P$
obtained by factoring through the diagonal $C \subset C \times _S C$; in
particular, this is also the identity when restricted to the diagonal in $C
\times ^{\PD}_S C$, so the difference of these two isomorphisms gives an
automorphism of $p_1^* P$ over $\Spec \O_C \oplus F_C^* \omega_{C/S}$, which
is the identity modulo $F_C^* \omega_{C/S}$. Because of the square-zero
structure on this sheaf, this corresponds to a section of $\InfAut(P)
\otimes F_C^* \omega_{C/S}$, where $\InfAut(P)$ is the sheaf of infinitesmal
automorphisms of $P$ over $C$. In the case of $P=\P(\E)$ for some vector
bundle $\E$, we have $\InfAut(P) = \End^0(\E)$. This section of $\End^0(\E)
\otimes F_C^* \omega_{C/S}$ is then our $p$-curvature associated to
$\nabla$. 

We remark that these definitions can be checked to agree with those of
\cite[A.1]{bo1} when the latter is defined, which is when everything (including the base scheme) is smooth over a field.

\begin{proof}[Proof of Proposition \ref{dormant}] We need only check the following two statements: first, that given any $\nabla$ on $\E$, the $p$-curvature of the induced connection on $\P(\E)$ is simply the traceless part of the usual $p$-curvature of $\nabla$ on $\E$; and second, that a connection $\nabla$ on $\E$ with $p$-trivial determinant has $p$-curvature taking values in $\End^0(\E) \otimes F_C^* \omega_{C/S}$.

We begin by checking that the new definition of $p$-curvature agrees with
the usual notion, when applied simply to $\E$. We can check this on local coordinates, so
say $t$ is a local coordinate for $C/S$, then the derivations are locally
generated by $\frac{d}{dt}$, so it suffices to check the assertion on this
derivation, which satisfies $(\frac{d}{dt})^p=0$. Thus, the usual definition
of $p$-curvature is simply $(\nabla_{\frac{d}{dt}})^p$. On the other hand,
the lifting of our connection to the $(p+1)$st order neighborhood in $C
\times^{\PD}_S C$ is given by the formula $\sum _{i=0} ^{p} (\nabla_{\frac{d}{dt}})^i \otimes (1 \otimes t - t \otimes 1)^{[i]}$, and when we mod out by $\I$, which is locally generated by $1 \otimes t - t \otimes t$, this kill all terms of positive order less than $p$, leaving only $1+(\nabla_{\frac{d}{dt}})^p \otimes (1 \otimes t - t \otimes 1)^{[p]}$. The tautological map in the opposite direction is given simply by $1$, so our automorphism of $p_1^* \E$ is given simply by $1+(1 \otimes t - t \otimes 1)^{[p]} \otimes (\nabla_{\frac{d}{dt}})^p$. But $(1 \otimes t - t \otimes 1)^{[p]}$ is precisely the square-zero element we used to obtain an infinitesmal automorphism of $\E$, so we find that the $p$-curvature from our definition is also $(\nabla_{\frac{d}{dt}})^p$, just as in the usual case.

One then easily checks our first assertion, since the functorial nature of the new definition of $p$-curvature implies that the section of $\End^0(\E) \otimes F_C^* \omega_{C/S}$ associated to the induced connection on $\P(\E)$ is simply the traceless part of the $p$-curvature of $\nabla$ on $\E$, since this is precisely the description of the natural map $\End(\E) \cong \InfAut(\E) \rightarrow \InfAut(\P(\E)) \cong \End^0(\E)$. The second assertion follows similarly, since the natural map $\End(\E) \cong \InfAut(\E) \rightarrow \InfAut(\det \E) \cong \O_C$ is precisely the trace map. 
\end{proof}

We now address our statements of Mochizuki's results.

\begin{proof}[Proof of Theorem \ref{moch1}] This is the $n=0$ case of \cite[Thm. II.2.8, p. 153]{mo3}.
\end{proof}

\begin{proof}[Proof of Theorem \ref{moch2}] We claim that it is enough to
see that the stack of dormant torally crys-stable bundles of level $\ell$ is
smooth over $\fF_p$ of relative dimension $5g-5+2r-2\ell$. Indeed, by
standard commutative algebra it then suffices to show that the fibers of the
map to $\overline{\M}_{g,r}$ have dimension at most $2g-2+r-2\ell$, which
follows from \cite[Thm. II.1.9, p. 132]{mo3}, after noting that the
condition of dormancy gives a closed sub-functor of the condition of 
nilpotency (defined to be nilpotency of the $p$-curvature of 
the connection), and that the 
condition on the characteristic insures that 
$(\overline{\D}^{\ell}_{g,r})'$ is in fact all of 
$\overline{\D}^{\ell}_{g,r}$. Note also the relevant statements on the map
$\overline{\D}^{\ell}_{g,r} \rightarrow \overline{\M}_{g,r}$ following
\cite[Def. I.3.7, p. 105]{mo3}. 

We have thus reduced down to the smoothness assertion. Over $\M_{g,r}$, this is the $j=0, i=\ell$ case of \cite[Cor. III.1.6, p. 176]{mo3}, but we will need the full statement, whose proof is more involved. Given $C_0, \E_0, \L_0, \nabla_0$ a dormant torally crys-stable bundle of level $\ell$ on $C_0$, the key is to consider the maps
\begin{equation}
\xymatrix{{\Def(C_0, \E_0, \nabla_0)} \ar[r] \ar[d] & {H^0(C_0, \End^0(\E_0)\otimes F^* \omega_{C_0^{(p)}})^{\nabla_0}} \\
{H^1(C_0, \E_0/\L_0 \otimes \L_0^{-1})} & {}}\end{equation}
where $\Def(C_0, \E_0, \nabla_0)$ denotes the space of logarithmic
first-order infinitesmal deformations of the triple, and $\omega$ the
logarithmic differentials. The horizontal arrow is the $p$-curvature map,
while the vertical arrow is the map giving the obstruction to deforming
$\L_0$ inside the given deformation of $\E_0$; see \cite[Prop. 1.7, p.
94]{mo3}.  We also remark that we consider deformations with a fixed
determinant, in the sense that if $\det \E_0 \cong F^* \sM_0$ for some
$\sM_0$ on $C_0^{(p)}$, we deform $\sM_0$ over the chosen deformation of
$C_0^{(p)}$, pull back under $F$, and require that the deformations of
$\E_0, \nabla_0$ have determinants equal to the resulting line bundle with
connection. The key point is if we have a logarithmic deformation $C_1$ of
$C_0$, the theory of crystals gives a natural identification between the
spaces $\Def_{C_1}(\E_0, \nabla_0)$ as $C_1$ is allowed to vary, so that we have
a splitting $\Def(C_0, \E_0, \nabla_0) = \Def(C_0) \times \Def(\E_0, \nabla_0)$,
where $\Def(\E_0, \nabla_0)$ is an abstract torsor over $H^1_{\DR}(C_0,
\End^0(\E_0) \overset{\nabla_0}{\rightarrow} \End^0(\E_0)\otimes \omega_{C_0})$. Moreover, under this splitting, the $p$-curvature map factors through $\Def(\E_0, \nabla_0)$.
One checks that as the element of $\Def(C_0)$ varies from a fixed choice, the obstruction map is the natural map $H^1(C_0, \tau_{C_0}) \rightarrow H^1(C_0, \E_0/\L_0 \otimes \L_0^{-1})$ induced by the Kodaira-Spencer map of $\nabla_0$ after tensoring by $\tau_{C_0}\otimes \L_0^{-1}$, where $\tau_{C_0}$ denotes the logarithmic tangent sheaf. Surjectivity of the obstruction map follows from the hypothesis that the Kodaira-Spencer map is non-zero, and hence has a torsion sheaf for its cokernel.

To complete the proof, it was pointed out to the author by Mochizuki that the arguments and results \cite[pp. 150-152]{mo3} do not in fact make use of the torally indigenous hypothesis, requiring only that the bundle be torally crys-stable. We then find that the $p$-curvature map is surjective by \cite[Lem. 2.7, p. 152]{mo3}. This surjectivity implies the desired smoothness, by standard deformation theory arguments over small extensions, as in the proof of smoothness of dormant torally indigenous bundles in the $n=0$ case of \cite[Thm. II.2.8, p. 153]{mo3}. Lastly, the dimension is the dimension of the simultaneous kernel of the $p$-curvature and obstruction maps. Again as in the proof of \cite[Thm. II.2.8, p. 153]{mo3}, the kernel of the $p$-curvature map on $\Def(\E_0, \nabla_0)$ has dimension $3g-3+r$,
and by the surjective of the obstruction map, we obtain the dimension of the total kernel by adding $\dim \Def(C_0)-\dim
H^1(C_0, \E_0/\L_0 \otimes \L_0^{-1})=3g-3+r-(g-1+2\ell)=2g-2+r-2\ell$, giving a total of $5g-5+2r-2\ell$, as desired.
\end{proof}

We conclude with two results stated entirely in terms of Mochizuki's language (in particular, they are over singular curves), but which we will need for our applications.

\begin{prop}\label{glue} A dormant torally indigenous bundle of radii $\{\rho_i\}$ on a nodal curve $C$ is equivalent to a collection of dormant torally indigenous bundles on each connected component $\tilde{C}_i$ of the normalization of the curve, where points on the $\tilde{C}_i$ lying above nodes are designated marked points, having radii $\{\rho_i\}$ at the points lying above marked points of $C$, and arbitrary radii at points lying above nodes, subject to the restriction that the radii of the two points lying above any given node must agree.
\end{prop}

\begin{proof} This follows immediately from \cite[\S I.4.4, p. 118]{mo3}, taking into account the fact that dormancy, being a condition of vanishing of $p$-curvature, will not be affected by the gluing.
\end{proof}

We remark that the following lemma, whose proof is trivial from the results 
of \cite{os10}, would be far more complicated if one attempted to directly 
apply Mochizuki's technique of considering totally degenerate curves, 
because such curves have no torally crys-stable bundles of intermediate level. In particular, by considering instead irreducible rational nodal curves, we avoid the machinery of PTCS bundles
\cite[\S II.1.6, p. 137]{mo3}.

\begin{lem}\label{nonempty} Let $C$ be an irreducible rational nodal curve of arithmetic genus $g$, with no marked points. Then for any $\ell>0$, there exists a dormant torally crys-stable bundle of level $\ell$ on $C$.
\end{lem}

\begin{proof} Let $\tilde{C}$ be the normalization of $C$, with $2g$ marked points lying above the nodes. We first consider the case that $\ell$ is an integer. Let $f$ be any function from $\tilde{C} \cong \P^1 \rightarrow \P^1$ of degree $g-\ell$, and unramified at the marked points. 
By
\cite[Thm. 6.7]{os10}, this gives a dormant torally crys-stable bundle of level $\ell$ on $\tilde{C}$ with radii $\frac{p-1}{2}$ at every marked point, and by the previous proposition we can glue to obtain a dormant torally crys-stable bundle of level $\ell$ on $C$, as desired. The case that $\ell$ is a half-integer is the same, except that (choosing an identification $\tilde{C}\cong \P^1$ such that $\infty$ is not a marked point), we need $f$ to have degree $p+g-\ell$, and be ramified to order at least $p$ at $\infty$. This is easily accomplished by choosing $f$ as before, with the additional hypotheses that $f(\infty)=0$ and $f$ is unramified at $\infty$, and then adding $x^p$.
\end{proof}

\section{Frobenius-Unstable Vector Bundles}

In this section, we apply Theorem \ref{moch2} to the study of the locus of
Frobenius-unstable vector bundles on a curve. In particular, we restrict
further to the case that $r=0$. We denote by $C^{(p)}$ the $p$-twist of $C$
over $S$, and Frobenius will always refer to the relative Frobenius map
$F:C \rightarrow C^{(p)}$. While it would be convenient to work with the entire Frobenius-unstable locus at once, it appears difficult to define this functorially, so we first define:

\begin{defn} Given a half-integer $\ell > 0$, a {\bf Frobenius-unstable vector bundle of level $\ell$} is a semi-stable vector bundle $\F$ of rank $2$ together with a line sub-bundle of $F^* \F$ with degree equal to $\ell+\frac{p}{2}\deg \F$.
\end{defn}

We then claim:

\begin{prop}\label{frob-unstab-equiv} Given a positive half-integer $\ell$, fix a line bundle $\L$ on $C^{(p)}$, with degree having the same parity as $2\ell$. Then dormant torally crys-stable bundles of level $\ell$ are equivalent to vector bundles on $C^{(p)}$ of determinant $\L$ which are Frobenius-unstable of level $\ell$, up to tensoring by (one of 
the $2^{2g}$) line bundles of order $2$.
\end{prop}

\begin{proof} 
First, given an $\F$ on $C^{(p)}$ which is semi-stable, but such that 
$F^* \F$ has a destabilizing sub-bundle of level $\L$, it is easy to see that
$(F^* \F, \Nc)$ gives a dormant torally crys-stable of level $\ell$, 
taking into account the comment following Proposition \ref{level}, and 
noting that the fact that $r=0$ renders the additional conditions
irrelevant.

Conversely, Proposition \ref{deproj} tells us that given a torally 
crys-stable bundle of level $\ell$, we can choose a representing pair 
$(\E, \nabla)$ such that $\E$ has determinant $F^* \L$, with $\nabla$ having 
determinant equal to the corresponding canonical connection. 
Next, given the dormancy condition, the Cartier isomorphism
(see \cite[Thm. 5.1]{ka1}) gives us an $\F$ on $C^{(p)}$ with determinant $\L$, such that $(\E, \nabla)= (F^* \F, \Nc)$. The destabilizing
sub-bundle of $\E$ in the definition of level of a torally crys-stable
bundle gives us the required destabilizing sub-bundle of $F^* \F$, but the
fact that the destabilizing sub-bundle is not horizontal 
for $\Nc$ implies by Lemma \ref{destab-unique} that $\F$ itself is stable. The only ambiguity in this was
in the choice of $\E$, since once $\nabla$ is required to have a particular determinant it is uniquely determined. But $\E$ is determined up to
tensoring by a line bundle, and for the determinant to remain unchanged, the
line bundle must be of order $2$. Since $p>2$, this is equivalent to
tensoring $\F$ by a line bundle of order $2$, as desired. 
\end{proof}

\begin{lem} The functor of Frobenius-unstable vector bundles of level $\ell$
is a locally closed sub-functor of the functor of stable vector bundles (of
rank 2 and the appropriate determinant).
\end{lem}

\begin{proof} By Lemma \ref{destab-unique}, Frobenius-unstable vector
bundles form a sub-functor of the functor of semi-stable vector bundles, and
we see by the same lemma that they must lie within the stable locus. We thus
need to check that this is a locally closed subfunctor; since both functors 
deal only with quasi-coherent sheaves which are coherent locally on the base,
and maps between such sheaves, they are clearly locally of finite type, so
we may restrict to schemes $T$ of finite type over $S$, and in particular 
Noetherian. Now, suppose we are given a 
vector bundle $\F$ of rank $2$ on $C^{(p)} \times_S T$; we want to show that
the locus of $T$ on which $F^* \F$ has a line sub-bundle, or equivalently,
locally free quotient of rank 1, of the required degree, is given
functorially as a locally closed subscheme of $T$. Let $Q/T$ be the
Quot scheme of quotients of $F^* \F$ of the required rank and 
degree, and $Q' \subset Q$ the open subscheme parametrizing locally free 
quotients; see \cite[Thm. 2.2.4, Lem. 2.1.7, Lem. 2.1.8]{h-l}, and note that
$Q$ is proper over $T$. We want to
show that $Q' \rightarrow T$ is an immersion. We already observed that it is
a monomorphism, so by \cite[Cor. I.2.13, p. 102]{mo3} it suffices to check
the valuative criterion for radimmersions, \cite[Thm. I.2.12, p. 101]{mo3}.
Because $Q$ is proper over $T$, we immediately see that it is enough to show
that any $k$-valued point of $Q$ whose image in $T$ is in the image of
$Q'$ must itself be in $Q'$. This is equivalent to the assertion that over a
field, if $F^* \F$ has a locally free quotient of rank 1 and the required
degree, then it has no other quotients of the same rank and degree. Finally,
we can conclude the desired statement because we are on a smooth curve, so
taking a saturation of the kernel of such a quotient would produce a 
destabilizing line sub-bundle, and we could then apply Lemma
\ref{destab-unique} once more.
\end{proof}

It thus follows via Theorem \ref{moch2} and general moduli space theory that 
we can conclude:

\begin{thm}\label{frob-unstab} Suppose that $p > 2g-2$, $S=\Spec k$, and $C$ is general.
Then the locus of Frobenius-unstable vector bundles of trivial determinant 
in the coarse moduli space $M_2(C^{(p)})$ of semi-stable vector bundles of 
rank two and trivial determinant on $C^{(p)}$ is non-empty of dimension $2g-4$. The locus of Frobenius-unstable vector bundles inside the coarse moduli space of semi-stable vector bundles of rank two and fixed odd determinant is non-empty of dimension $2g-3$.
\end{thm}

\begin{proof} We may work level by level, since the Frobenius-unstable locus
in the moduli space is clearly the image of the disjoint union of the loci
for each level. Noting that the argument of 
Proposition \ref{frob-unstab-equiv} used no hypotheses on the base to go 
from a
Frobenius-unstable bundle to a dormant torally crys-stable bundle, we get a morphism 
from the functor of Frobenius-unstable bundles of level $\ell$ to the algebraic space of dormant torally
crys-stable bundles of level $\ell$. Moreover, Proposition
\ref{frob-unstab-equiv} implies that this morphism is an isomorphism on
strictly Henselian schemes over fixed geometric points, so we conclude that
the morphism is formally \'etale. 

Since the stable locus of coarse moduli space $M_2(C^{(p)})$ represents the 
\'etale sheafification of the corresponding functor, the locus of Frobenius-unstable bundles of level $\ell$ represents the sheafification of the functor of Frobenius-unstable bundles of level $\ell$. In particular any map to an algebraic space factors through the locus in question. Furthermore, the map from the moduli functor of semi-stable vector bundles to $M_2(C^{(p)})$ is bijective on $k$-valued points and formally \'etale, so we conclude that there is a map from the locus of Frobenius-unstable bundles of level $\ell$ in $M_2(C^{(p)})$ to the algebraic space of dormant torally crys-stable bundles of level $\ell$, and this map is formally \'etale, hence \'etale. We thus conclude the dimension of the locus of Frobenius-unstable bundles of level $\ell$ from Theorem \ref{moch2}. We may finish the proof by considering the minimal cases $\ell=\frac{1}{2}, 1$ as appropriate, as long as we know that the spaces in question are non-empty. By the flatness assertion of Theorem \ref{moch2}, it suffices to show that the locus of torally crys-stable bundles of level $\ell$ is non-empty for all $\ell>0$ on any single curve $C_0$, which is Lemma \ref{nonempty}.
\end{proof}

\margh{check odd det case for $g=2$; diff for $p=2$?}

We also mention that in the only case for which the number of Frobenius-unstable bundles is finite, Mochizuki's results also give the number.

\begin{thm} The number of Frobenius-unstable vector bundles of rank $2$ and
trivial determinant on a general curve $C$ of genus $2$ is given by
$\frac{2(p^3-p)}{3}$. Further, they all have no non-trivial deformations, in
the sense that for a Frobenius-unstable bundle $\F$, no non-trivial 
first-order deformation of $\F$ induces the trivial deformation of $F^* \F$.
\end{thm}

\begin{proof} In this case, the only possible positive integral level is
$\ell=1$. By Theorem \ref{frob-unstab-equiv}, the number Frobenius unstable
bundles is $2^{2g}=16$ times the number of dormant torally indigenous
bundles on $C$. By Theorem \ref{moch1}, there are no non-trivial deformations
for a general $C$, and it suffices to consider the case that $C$ is totally degenerate; we choose the totally degenerate curve obtained by gluing two copies of $\P^1$ to each other three times. By Proposition \ref{glue}, and the fact
\cite[p. 206]{mo3} that a torally indigenous bundle on $\P^1$ with three marked points is determined by its radii, we find that this number is the same as the number of dormant torally indigenous bundles on $\P^1$ with three marked points, which is $\frac{p^3-p}{24}$ by
\cite[Cor. V.3.7, p. 267]{mo3}.
\end{proof}

\begin{rem} One natural question which is not addressed by these techniques, due to the necessity of working level by level, is the relationship between the strata of different levels. Particularly if one wanted to conclude, for instance, that the Frobenius-unstable locus is pure of the given dimension, one would want to show that every higher-level bundle is a specialization of one with lower (but still positive) level. However, the deformation theory raised by such a question seems substantially more delicate than the constant-level case.
\end{rem}

\section{Rational Functions with Prescribed Ramification}\label{s-maps}

For the desired applications to self-maps of $\P^1$, the statements of
equivalences with dormant torally indigenous bundles are slightly more
complicated. We specialize to the case $g=0$, and show:  

\begin{thm}\label{map-equiv} Fix radii $\{\rho_i\}$ for each of the $r$
marked points on $C=\P^1$, with sign mod $p$ chosen so that
$0<\rho_i<\frac{p}{2}$. Then there is a natural map $\vp$ from the set of
self-maps of $\P^1$ ramified to order $p-2\rho_i$ at the marked points and
unramified elsewhere, to the set of dormant torally indigenous bundles on
$\P^1$ of radii $\{\rho_i\}$. We have further:
\begin{ilist}
\itm The map $\vp$ is injective after passing to equivalence classes of maps related by post-composition by automorphisms of the image.
\itm If the marked points are general, the map $\vp$ is surjective.
\itm If the marked points are general, there is a bijective correspondence between self-maps of $\P^1$ as above, and self-maps of $\P^1$ for which any even number of the ramification indices $p-2\rho_i$ have been replaced by $2\rho_i$.
\end{ilist}
\end{thm}

\begin{proof} Since we are on $\P^1$, there need not be any ambiguity in our
choices of vector bundles with connection: for $r$ even, we can choose $\E =
\O(-\frac{r}{2}+1) \oplus \O(\frac{r}{2}-1)$, and for $r$ odd, we can choose 
$\E = \O(\frac{p-r}{2}+1) \oplus \O(\frac{r+p}{2}-1)$, and then $\nabla$ will 
be uniquely determined as there is only one connection with vanishing 
$p$-curvature on any line bundle. We then apply \cite[Thm. 1.1]{os10} to 
conclude the existence of $\vp$ as well as assertions (i) and
(ii). Assertion (iii) is in fact unrelated to indigenous bundles; this is
simply \cite[Lem. 5.2]{os7}.
\end{proof}

One conclusion is:

\begin{cor} The number of self-maps of $\P^1$ ramified to order less than $p$ at $r$ general marked points $P_i$ and unramified elsewhere, counted modulo automorphism of the image $\P^1$, is $2^{r-1}$ times the number of dormant torally indigenous bundles on $\P^1$ with the same marked points.
\end{cor}

\begin{proof} In light of the previous theorem, it suffices to check that
$2^{r-1}$ gives the number of ways of replacing an even number of the
$p-2\rho_i$ by $2 \rho_i$, which is to say, the number of even subsets of
$\{1, \dots, r\}$. But this is simple enough: for instance, one sees by
expanding $(1-1)^r$ with the binomial theorem that the number of even subsets is the same as the number of odd ones.
\end{proof}

More importantly, as an immediate corollary of (i) of our theorem together with the finiteness assertion of Theorem \ref{moch1}, we can conclude:

\begin{thm}\label{finite} Fix distinct points $P_1, \dots, P_r$ on $\P^1$, and odd integers $e_1, \dots, e_r$ less than $p$. Then the number of maps from $\P^1$ to itself which are ramified to order $e_i$ at the $P_i$ and unramified elsewhere, when counted modulo automorphisms of the image $\P^1$, is finite.
\end{thm}

\begin{rem} The parity hypothesis on the $e_i$ in this theorem is not
expected to be necessary; indeed, if based on the desired parity of
ramification indices one varies the determinant of the vector bundle with
connection corresponding to the given torally indigenous bundle, as in
\cite[p. 206]{mo3}, one should obtain the injectivity of Theorem 
\ref{map-equiv} (i) without parity hypotheses. In contrast, the hypothesis that
all $e_i$ are less than $p$ is unquestionably necessary, as demonstrated by
families such as $x^{p+2}+tx^p+x$. We give a direct proof of this finiteness
result when the $P_i$ are general (and without assumptions on the parity or
size of the $e_i$, as long as $p \nmid e_i$) in
\cite{os7}, but a direct proof for arbitrary distinct $P_i$ remains elusive. The proof in the context of dormant torally indigenous bundles is carried out by first enlarging to the category of nilpotent torally indigenous bundles, and proving a finiteness result there; it is this construction which is lacking in the context of self-maps of $\P^1$.
\end{rem}

This finiteness result leads to strong non-existence statements for certain
tame branched covers of $\P^1$ by $\P^1$; see \cite{os12}.

We conclude with a discussion of the relevance of dormant torally indigenous bundles to counting self-maps of $\P^1$ with prescribed ramification. Specifically, Proposition \ref{glue} allows us to obtain an alternate proof of the main theorem of \cite{os7}. 

\begin{thm} Given $d, n$ and $e_1, \dots, e_n$ with $e_i<p$ and $e_i \leq d$ for all $i$ and $2d-2 = \sum_i (e_i-1)$, the number $N_{\gen}(\{e_i\}_i)$ of separable self-maps of $\P^1$ of degree $d$ and ramified to order $e_i$ at $P_i$, modulo automorphism of the image space, is given by

\begin{equation}N_{\gen}(\{e_i\}_i) = 
\!\!\!\!\!\!\!\!\!\!\!\!\!\!\!\!\!\!\!\!\!\!
\sum
_{\scriptsize \begin{matrix}d-e_{n-1}+1 \\ d-e_{n}+1\end{matrix} \leq d' \leq
\begin{matrix}d \\ p+d-e_{n-1}-e_n\end{matrix}} 
\!\!\!\!\!\!\!\!\!\!\!\!\!\!\!\!\!\!\!\!\!\!
N_{\gen}(\{e_i\}_{i \leq n-2}, e), \text{ with } e=2d'-2d+e_{n-1}+e_n-1
\end{equation}

\begin{equation}N_{\gen}(e_1, e_2, e_3)= \begin{cases}1 & p>d \\ 0 &
\text{otherwise} \end{cases}
\end{equation}

Equivalently, $N_{\gen}(\{e_i\}_i)$ is given as the number of ways of inserting $\lceil \frac{n}{3} \rceil$ indices $e'_j$ into the sequence $e_1, \dots, e_n$, starting after $e_2$ and alternating thereafter, such that any consecutive triple $e, e', e''$ of the newly obtained sequence satisfies the triangle inequality, and we further have that $e+e'+e''$ is odd and less than $2p+1$.
\end{thm}

\begin{proof} The two statements are easily seen to be equivalent by induction on $n$. By 
\cite[Cor. 5.3]{os7}, it suffices to handle the case where all $e_i$ are odd. By Theorem \ref{map-equiv}, we see that $N_{\gen}(\{e_i\}_i)$ then counts the number of dormant torally indigenous bundles of radii $\{\frac{p-e_i}{2}\}_i$. By Theorem \ref{moch1}, it may then be computed on a totally degenerate curve; we choose a chain with $P_1, P_2$ lying on the first component, $P_3$ on the second, and so forth, until $P_{n-1}$ and $P_n$ lie on the last component. By Proposition \ref{glue}, and then another application of Theorem \ref{map-equiv}, this may then be described precisely as claimed in our second formulation, using also our description of the situation for $n=3$.
\end{proof}

Although this proof is far more circuitous in translating to
vector bundles with connection, invoking Mochizuki's theory, and then
translating back to maps, both chronologically and in spirit it preceded the
ultimately self-contained proof of \cite{os7}. Indeed, although a self-contained argument was the original intent, after
some study of the situation there remained the stumbling block of
controlling degeneration from separable to inseparable maps, ultimately
addressed by \cite[Thm. 6.1]{os7}. However, we know from 
\cite[Thm. 6.7]{os10} that for the special configurations of marked points where the map $\vp$ of Theorem \ref{map-equiv} fails to be surjective, the dormant torally indigenous bundles which are not in the image in fact correspond to families of self-maps of $\P^1$ with the prescribed ramification at the marked points, together with an additional ramification point of order greater than $p$ at an additional point. By 
\cite[Prop. 5.4]{os7} such families are known to exist only for special configurations (and indeed, this is how the surjectivity of $\vp$ is shown for general configurations of marked points). On the other hand, separable maps can certainly degenerate to inseparable maps as the ramification points are allowed to move. Now, Mochizuki's finite flatness theorem implies that dormant torally indigenous bundles can neither appear nor disappear under specialization, so the natural conclusion is that separable maps can degenerate to inseparable maps only when the ramification points move into special configurations allowing the existence of families of maps with exactly one ramification index greater than $p$. With this realization, it was then possible to write down examples of degenerations of connections which explicitly demonstrated the phenomenon, and finally to discover the construction which led to a self-contained proof of this statement
\cite[Thm. 6.1]{os7} and hence the formula for self-maps of $\P^1$ with prescribed ramification.

\bibliographystyle{hamsplain}
\bibliography{hgen}

\end{document}

%% file: headerstemplate.tex
\usepackage{amssymb,euscript,amsmath, mathrsfs}
\usepackage[dvips]{graphicx}
\usepackage[dvips]{color}

\newcounter{ENUM}
\newcommand{\itm}{\item}
\newenvironment{ilist}{\renewcommand{\theENUM}{\roman{ENUM}}\renewcommand{\itm}{\addtocounter{ENUM}{1}\item[(\theENUM)]}\begin{itemize}\setcounter{ENUM}{0}}{\end{itemize}}

\input xy
\xyoption{all}
\CompileMatrices

\def\P{{\mathbb P}}

\def\fF{{\mathbb F}}
\def\E{{\mathscr E}}
\def\F{{\mathscr F}}
\def\sM{{\mathscr M}}

\def\D{{\mathcal D}}

\def\L{{\mathcal L}}
\def\O{{\mathcal O}}
\def\I{{\mathcal I}}
\def\J{{\mathcal J}}

\def\M{{\mathcal M}}

\def\vp{\varphi}

\def\Nc{\nabla ^{\text{can}}}

\def\et{\text{\'et}}

\def\gen{\operatorname{gen}}

\def\Hom{\operatorname{Hom}}

\def\End{\operatorname{End}}

\def\Ad{\operatorname{Ad}}

\def\Tr{\operatorname{Tr}}
\def\Spec{\operatorname{Spec}}

\def\Res{\operatorname{Res}}

\def\PConn{\operatorname{PConn}}

\def\Def{\operatorname{Def}}
\def\PD{\operatorname{PD}}
\def\DR{\operatorname{DR}}
\def\InfAut{\operatorname{InfAut}}

\newcommand{\margh}[1]{}

\newtheorem{thm}{Theorem}[section]
\newtheorem{prop}[thm]{Proposition}
\newtheorem{lem}[thm]{Lemma}
\newtheorem{cor}[thm]{Corollary}

\theoremstyle{definition}
\newtheorem{defn}[thm]{Definition}

\theoremstyle{remark}

\newtheorem{rem}[thm]{Remark}

\numberwithin{equation}{section}

%% file: crys-lex.bbl
\newcommand{\noopsort}[1]{} \newcommand{\printfirst}[2]{#1}
  \newcommand{\singleletter}[1]{#1} \newcommand{\switchargs}[2]{#2#1}
\providecommand{\bysame}{\leavevmode\hbox to3em{\hrulefill}\thinspace}
\begin{thebibliography}{10}

\bibitem{b-o}
Pierre Bertholot and Arthur Ogus, \emph{Notes on crystalline cohomology},
  Princeton University Press, 1978.

\bibitem{bo1}
Jean-Benoit Bost, \emph{Algebraic leaves of algebraic foliations over number
  fields}, Institut Des Hautes Etudes Scientifiques Publications Mathematiques
  (2001), no.~93, 161--221.

\bibitem{h-l}
Daniel Huybrechts and Manfred Lehn, \emph{The geometry of moduli spaces of
  sheaves}, Max-Planck-Institut fur Mathematik, 1997.

\bibitem{j-r-x-y}
Kirti Joshi, S.~Ramanan, Eugene~Z. Xia, and Jiu-Kang Yu, \emph{On vector
  bundles destabilized by {F}robenius pull-back}, preprint.

\bibitem{ka4}
Kazuya Kato, \emph{Logarithmic structures of {F}ontaine-{I}llusie}, Algebraic
  Analysis, Geometry, and Number Theory, John Hopkins Univ. Press, 1989,
  pp.~191--224.

\bibitem{ka1}
Nicholas~M. Katz, \emph{Nilpotent connections and the monodromy theorem:
  Applications of a result of {Turrittin}}, Inst. Hautes Etudes Sci. Publ.
  Math. \textbf{39} (1970), 175--232.

\bibitem{l-p3}
Herbert Lange and Christian Pauly, \emph{On {F}robenius-destabilized rank-2
  vector bundles over curves}, \mbox{arXiv:math.AG/0309456}.

\bibitem{l-p}
Yves Laszlo and Christian Pauly, \emph{The action of the {F}robenius maps on
  rank 2 vector bundles in characteristic 2}, Journal of Algebraic Geometry
  \textbf{11} (2002), no.~2, 129--143.

\bibitem{mo3}
Shinichi Mochizuki, \emph{Foundations of $p$-adic {T}eichm\"uller theory},
  American Mathematical Society, 1999.

\bibitem{og2}
Arthur Ogus, \emph{F-crystals, {G}riffiths transversality, and the {H}odge
  decomposition}, Ast\'erisque \textbf{221} (1994).

\bibitem{os11}
B.~Osserman, \emph{{F}robenius-unstable bundles and $p$-curvature},
  \mbox{arXiv:math.AG/0409266}.

\bibitem{os12}
\bysame, \emph{Linear series and existence of branched covers}, in preparation.

\bibitem{os10}
\bysame, \emph{Logarithmic connections with vanishing $p$-curvature},
  \mbox{arXiv:math.AG/0409145}.

\bibitem{os7}
\bysame, \emph{Rational functions with given ramification in characteristic
  $p$}, \mbox{arXiv:math.AG/0407445}.

\end{thebibliography}
